\documentclass[a4paper,oneside,reqno,12pt]{amsart}

\usepackage{cmap}
\usepackage[T1]{fontenc}
\usepackage[utf8]{inputenc}
\usepackage[english]{babel}
\usepackage{csquotes}
\usepackage{lmodern}
\usepackage{amsmath,amssymb,amsfonts,mathtools}
\usepackage[
    a4paper,
    hmargin=2.2cm,
    vmargin=2.8cm
]{geometry}
\usepackage{color, xcolor}

\definecolor{myblue}{HTML}{0072B2}
\definecolor{myorange}{HTML}{E69F00}
\definecolor{mygreen}{HTML}{009E73}
\definecolor{myred}{HTML}{D55E00}
\definecolor{mypurple}{HTML}{CC79A7}
\definecolor{myskyblue}{HTML}{56B4E9}

\usepackage{hyperref}
\hypersetup{
    colorlinks = true,
    urlcolor   = myblue,
    linkcolor  = myred,
    citecolor  = mygreen
}

\newtheorem{theorem}{Theorem}
\newtheorem{lemma}[theorem]{Lemma}

\DeclareMathOperator{\inter}{int}
\newcommand{\R}{\mathbb{R}}

\title{A sharp plank theorem for cube annuli}

\author[E. Bakaev]{Egor Bakaev}
\address{Egor Bakaev,
Department of Computer Science, University of Copenhagen}

\author[A. Polyanskii]{Alexander~Polyanskii}
\address{Alexander Polyanskii,
Department of Mathematics, Emory University,
400 Dowman Dr, Atlanta, GA 30322, USA}
\email{\href{mailto:apolian@emory.edu}{apolian@emory.edu}}
\urladdr{\url{https://polyanskii.com}}


\begin{document}

\begin{abstract}
We prove a conjecture of Andr\'as Bezdek on plank coverings of square annuli. 
A~plank in a Euclidean space is the closed region between
two distinct parallel hyperplanes, and its width is the distance
between them. Let $Q$ be a unit square and let $H$ be a square of side
length $1/2$ whose sides are parallel to those of $Q$.
Then every finite family of planks covering $Q\setminus H$ has total
width at least~$1$. More generally, if an axis-parallel 
box with side lengths $h_1,\ldots,h_d$ is removed from the unit cube
in $\R^d$ and $\sum_jh_j = d-1$, then every finite plank covering of the remaining set  has total width at least one.
The bounds are sharp.
\end{abstract}

\maketitle

\section{Introduction}

A \emph{plank} in a Euclidean space is the closed region between
two distinct parallel hyperplanes, and its width is the distance
between them. The corresponding open region is called an
\emph{open plank}.
For a convex body $K$, its width $w(K)$ is the infimum
of the widths of planks containing~$K$. 

The classical plank problem, attributed to
Tarski~\cite{tarski1932uwagi}, asks whether $K$ can be covered by planks
whose total width is smaller than $w(K)$.
Bang~\cite{bang1950covering,bang1951solution} proved that this is impossible:
if finitely many planks cover $K$, then their total width is at least $w(K)$.

Andr\'as Bezdek~\cite{bezdek2003covering} asked whether the same lower
bound holds when a small homothetic copy of $K$ is removed from its
interior; see also Conjecture 7.3 in K\'aroly~Bezdek's
survey~\cite{bezdek2013tarski}.
The answer depends on $K$. White and Wisewell~\cite{white2007covering}
characterized the convex polygons for which every sufficiently small
homothetic hole leaves the infimum of the total covering width unchanged;
see also further developments in~\cite{ambrus2025covering}.
Bakaev and Yehudayoff~\cite{bakaev2026volumestabilityplankcovers}
proved Bezdek's conjecture for the concentric disk annulus as a
consequence of a more general theorem valid for all centrally symmetric
convex bodies. For the square, the remaining question is how large
the hole can be.

A \emph{box} in $\R^d$ is a Cartesian product of $d$ closed intervals.
All rectangles, squares, and cubes considered below are boxes
in this sense.

Let $Q$ be a unit square in the plane, and let $H$ be a square of side
length $h$.
A.~Bezdek conjectured that every plank covering of $Q\setminus\inter H$
has total width at least one when $h\le1/2$; see Remark 2
in~\cite{bezdek2003covering} and Conjecture 4 from Section 3.4
in~\cite{brass2005research}.
He proved the conjecture for $h\le1-1/\sqrt2$~\cite{bezdek2003covering}.
We settle the full conjecture.

\begin{theorem}\label{thm}
Let $Q$ be a unit square in $\R^2$ and let $H$ be a square
with side length $1/2$.
If finitely many planks cover $Q\setminus\inter H$, then their total width is at least one.
\end{theorem}

We also prove an extension of this result for an arbitrary dimension.
Smurov, Bogataya, and Bogatyi~\cite{smurov2010covering} proved the
corresponding result for a unit cube in $\R^d$ with a homothetic
cubical hole of side length $h\le1-1/\sqrt d$.
Recently, for $d=2,3$, the range was extended to
$h\le1-1/\sqrt d+\varepsilon_d$ for some particular $\varepsilon_d>0$;
see~\cite{ambrus2026plankcoveringscubeannuli}.
We prove the sharp bound $h\le1-1/d$, and more generally the following
statement for rectangular holes; it is easy to see that the
bound is sharp.

\begin{theorem}\label{thm:rectangular-hole}
Let $d\ge2$, let $Q$ be a unit cube in $\R^d$, and let
$H$ be a box whose side lengths $h_1,\ldots,h_d$ satisfy
\(
 \sum_{j=1}^d h_j = d-1.
\)
If finitely many planks cover $Q\setminus\inter H$, then their total width is at least one.
\end{theorem}

Our proof uses Bang's lemma~\cite{bang1951solution}, whose proof we include in the next section for completeness.
The idea of our proof is the following.
First, we make a suitable affine change of coordinates, apply
Bang's lemma to the transformed planks and find a point avoiding them, and then map the resulting
point back to the original coordinates. Here we vary only the coordinate scales. We choose them so that all signed sums in Bang's lemma,
mapped back to the original coordinates, fit inside the cube,
leaving enough room in one coordinate direction to add
an open plank containing the interior of the hole. The last part is the only technical challenge in the proof.

The same idea underlies Ball's work~\cite{ball1991plank}
and the short proof of his theorem by Akopyan and
Polyanskii~\cite{akopyan2026shortproofballsplank}.
In those proofs, the optimal affine transformation is chosen
according to a different criterion.

\section{Bang's lemma and the geometric reduction}\label{sec:geometry}

Let $P_1,\dots, P_m$ be open planks in $\mathbb R^d$. Write them as
\[
 P_i=\{z\in\R^d:|\langle n_i,z\rangle-s_i|< w_i/2\},
\]
where $n_i$ is a unit normal,
$w_i>0$ is the width, and $s_i\in\R$ determines the middle
hyperplane $\langle n_i,z\rangle=s_i$.
We use the following form of Bang's lemma.

\begin{lemma}[Bang's lemma~\cite{bang1951solution}]\label{bang}
There is a sign vector $\varepsilon=(\varepsilon_1,\dots, \varepsilon_m)\in\{-1,1\}^m$ such that
\[
 z_\varepsilon=\frac12\sum_{i=1}^m\varepsilon_iw_in_i
 \quad\text{and}\quad
 |\langle n_i,z_\varepsilon\rangle-s_i|\ge w_i/2 \qquad\text{for all }i.
\]
\end{lemma}
\begin{proof}
Choose $\varepsilon\in\{-1,1\}^m$ maximizing
\[
 E(\varepsilon):=\|z_\varepsilon\|^2
                  -\sum_{i=1}^m\varepsilon_iw_is_i.
\]
If
$\varepsilon^{(i)}$ is obtained by reversing the $i$th sign, then
$\varepsilon_i(\langle n_i,z_\varepsilon\rangle-s_i)\ge w_i/2$ for
every $i$, since
\[
 0\le E(\varepsilon)-E(\varepsilon^{(i)})
 =2w_i\bigl(\varepsilon_i(\langle n_i,z_\varepsilon\rangle-s_i)-w_i/2\bigr).
\qedhere
\]
\end{proof}
We assume that $Q$ is the unit cube centered at the origin.
For each choice of parameter vector $r=(r_1,\ldots,r_d)$,
we rescale the coordinate directions and apply Bang's lemma
to the transformed planks. Mapping the signed sums back to
the original coordinates gives a finite set of points, at
least one of which lies outside every plank.
We call the smallest box containing these points
\emph{Bang's box} for $r$.

Our first goal is to choose $r$ so that Bang's box is contained
in $Q$. If the total width of the planks is at most one,
this already holds for $r=(1,\ldots,1)$, since each side length
of Bang's box is at most the total width.
We also need enough room to add an open plank containing $\inter H$. 
If this plank has a coordinate normal, adding it
enlarges Bang's box only in that coordinate direction, by
the width of the added plank. We therefore seek parameters
for which the enlarged box is still contained in $Q$.
Bang's lemma then gives a point in $Q$ outside all 
the planks, including the added one, and hence outside $\inter H$.

We now describe Bang's box explicitly.
For a positive parameter vector $r=(r_1,\ldots,r_d)$ set 
\[
 A_r:=\operatorname{diag}(\sqrt{r_1},\ldots,\sqrt{r_d}).
\]

Consider the preimages of the planks
under the linear map
$z\longmapsto A_rz$, that is, the planks 
\begin{align*}
    A_r^{-1}P_i= & \{A_r^{-1}z\in \mathbb R^d: |\langle n_i, z\rangle -s_i|< w_i/2\}\\
    = & \{t\in \mathbb R^d:|\langle A_r n_i,t\rangle -s_i|< w_i/2\}\\
    = & \big\{t\in \mathbb R^d: \Big|\big\langle \frac{A_r n_i}{\|A_r n_i\|}, t\big\rangle - \frac{s_i}{\|A_r n_i\|}\Big| < \frac{w_i}{2\|A_rn_i\|}\Big\}.
\end{align*}
Apply Lemma~\ref{bang} to the planks $A_r^{-1}P_i$, whose widths
and unit normals are $w_i/\|A_rn_i\|$ and
$A_rn_i/\|A_rn_i\|$, respectively.
We obtain signs $\varepsilon_i\in\{-1,1\}$ such that the point
\[
 \sum_{i=1}^m\varepsilon_i
 \frac{w_i}{2\|A_rn_i\|}
 \frac{A_rn_i}{\|A_rn_i\|}
\]
lies outside every $A_r^{-1}P_i$.
Applying $A_r$ to this point, we conclude that one of the points
\[
 \sum_{i=1}^m\varepsilon_i v_i(r),
 \qquad \varepsilon_i\in\{-1,1\},
\]
lies outside every $P_i$, where
\begin{equation}\label{eq:signed-sums}
 v_i(r):=A_r\left(
       \frac{w_i}{2\|A_rn_i\|}
       \frac{A_rn_i}{\|A_rn_i\|}\right)
       =\frac{w_iA_r^2n_i}{2\|A_rn_i\|^2}.
\end{equation}
Note that if $n_i=\pm e_k$ for a standard basis vector $e_k$, then
$v_i(r)=\frac{w_i}{2}n_i$, independently of $r$.
The corresponding Bang's box is centered at the origin and has side lengths
\begin{equation*}\label{eq:coordinate-spans}
 X_k(r):=2\sum_{i=1}^m|(v_i(r))_k|
        =\sum_{i=1}^m\frac{w_i r_k |n_{ik}|}{\|A_r n_i\|^2}
        \qquad(1\le k\le d).
\end{equation*}

The following lemma gives the conditions under which we can
add an open plank containing $\inter H$ while keeping Bang's box inside a unit cube $Q$.

\begin{lemma}\label{lem:uncovered-point}
Let $H$ be a  box in $\R^d$
with side lengths $h_1,\ldots,h_d$, and let $P_1,\ldots,P_m$ be open planks.
If there is a positive parameter vector $r=(r_1,\ldots,r_d)$
and an index $j$ such that
\begin{equation*}\label{eq:uncovered-criterion}
 X_j(r)+h_j\leq 1 
 \quad\text{and}\quad
 X_k(r)\leq 1\quad\text{ for all }k\ne j,
\end{equation*}
then the planks do not cover $Q\setminus\inter H$.
\end{lemma}

\begin{proof}
If $h_j=0$ for some $j$, then $\inter H=\varnothing$ and Bang's box is
contained in $Q$, so the conclusion follows.
Assume that $h_j>0$ for all $j$.
Let $P_0$ be an open plank of width $h_j$ with normal $e_j$
containing $\inter H$.
By the observation written after~\eqref{eq:signed-sums}, its Bang vector is
$v_0(r)=\frac{h_j}{2}e_j$.
The signed sums for the family $P_0,P_1,\ldots,P_m$ are
\[
 \sum_{i=0}^m\varepsilon_i v_i(r),
 \qquad \varepsilon_i\in\{-1,1\}.
\]
They lie in the box centered at the origin with side lengths
$X_j(r)+h_j\le1$ in direction $j$ and $X_k(r)\le1$
in the other directions. Thus all these sums belong to $Q$.

Apply Lemma~\ref{bang} to the preimages of $P_0,\ldots,P_m$.
The image of the resulting point lies outside every $P_i$,
including $P_0$. Since $\inter H\subset P_0$, this point
belongs to $Q\setminus\inter H$.
\end{proof}

\section{The planar case}\label{sec:planar}

\begin{proof}[Proof of Theorem~\ref{thm}]

Suppose that closed planks of total width less than one cover
$Q\setminus\inter H$. Since $Q$ is bounded, we may slightly
rotate and enlarge them to obtain open planks whose unit normals
have no zero coordinates, while preserving the covering and
keeping the total width less than one.
Enlarging them further, we may assume that
$\sum\nolimits_i w_i=1$.
We retain the notation $P_i,w_i,n_i$ for these open planks.

Set $a_i:=|n_{i1}|$ and $b_i:=|n_{i2}|$.
Then $a_i,b_i>0$ and $a_i^2+b_i^2=1$.
Interchanging the coordinate axes if necessary, assume that
\begin{equation}
\label{equation: assumption}
 \sum\nolimits_i w_i\frac{a_i}{a_i+b_i}\le\frac 12 \sum\nolimits_i w_i=\frac12.
\end{equation}
The height $X_2(1,\lambda)$ is continuous and strictly
increasing in $\lambda$.
Moreover,
\[
\begin{aligned}
 X_2(1,1)
 &=\sum\nolimits_i\frac{w_i b_i}{a_i^2+b_i^2}
  =\sum\nolimits_i w_i b_i<\sum\nolimits_i w_i,\\
 \lim_{\lambda\to\infty}X_2(1,\lambda)
 &=\lim_{\lambda\to\infty}
   \sum\nolimits_i\frac{w_i\lambda b_i}{a_i^2+\lambda b_i^2}
 =\sum\nolimits_i\frac{w_i}{b_i}>\sum\nolimits_i w_i.
\end{aligned}
\]
Choose $\lambda>1$ so that
\begin{equation*}\label{eq:planar-balance}
 X_2(1,\lambda)
 =\sum\nolimits_i w_i\frac{\lambda b_i}{a_i^2+\lambda b_i^2}=\sum\nolimits_i w_i=1.
\end{equation*}
Apply Lemma~\ref{lem:two-variable} to each pair $(a_i,b_i)$, sum with weights $w_i$, and obtain that
\[
(\lambda-1)X_1(1,\lambda) + X_2(1,\lambda)
 \le \sum\nolimits_i w_i + (\lambda-1)\sum\nolimits_i w_i\frac{a_i}{a_i+b_i}.
\]
Using~the preceding equality and $\lambda>1$, we conclude
\[
X_1(1,\lambda) \leq \sum\nolimits_i w_i\frac{a_i}{a_i+b_i} \overset{\eqref{equation: assumption}}{\le}\frac12.
\]
So, we can apply Lemma~\ref{lem:uncovered-point} with $j=1$ and obtain a contradiction with the covering.
\end{proof}

\section{
A two-variable estimate
}

\begin{lemma}\label{lem:two-variable}
Let $a,b>0$ with $a^2+b^2\ge1$, and let $\lambda>1$.
Then
\begin{equation}\label{eq:two-variable-bound}
 \frac{(\lambda-1)a}{a^2+\lambda b^2}+\frac{\lambda b}{a^2+\lambda b^2}
 \le 1+(\lambda-1)\frac{a}{a+b}.
\end{equation}
\end{lemma}

\begin{proof}
The function
\[
 F(\lambda):=\frac{(\lambda-1)a+\lambda b}{a^2+\lambda b^2}
\]
has decreasing derivative
$F'(\lambda)= \frac{a(a^2+ab+b^2)}{(a^2+\lambda b^2)^2},$
so it is concave. Its tangent at $\lambda_0=1+a/b$ is~$
\frac{\lambda a+b}{a^2+ab+b^2}.
$
Also,
$a+b \leq ((a+b)^2+1)/2 \leq a^2+ab+b^2$.
Consequently,
\[
 \frac{(\lambda-1)a+\lambda b}{a^2+\lambda b^2}
 \le 
 \frac{\lambda a+b}{a^2+ab+b^2}    
 \le \frac{\lambda a+b}{a+b}
 =1+(\lambda-1)\frac{a}{a+b}.
 \qedhere
\]
\end{proof}

\section{Rectangular holes in arbitrary dimension}\label{sec:balancing}

Write $c_{ik}:=|n_{ik}|$.
The width of the unit cube in direction $n_i$ is $\sum_k c_{ik}$,
with $c_{ik}$ contributed by its $k$th coordinate edge.
Split each plank width in these proportions, defining
\begin{equation*}\label{eq:reference-widths}
 R_k:=\sum_i\frac{ c_{ik}}{\sum_{\ell=1}^d c_{i\ell}} w_i,
 \qquad \sum_{k=1}^d R_k=\sum_iw_i.
\end{equation*}

\begin{proof}[Proof of Theorem~\ref{thm:rectangular-hole}]
As in the proof of Theorem~\ref{thm}, a hypothetical covering
of total width less than one gives open planks with
$\sum_iw_i=1$ and $c_{ik}>0$ for every $i,k$,
which still cover $Q\setminus\inter H$.
Since
$
 \sum_jR_j=1=\sum_j(1-h_j),
$
some coordinate $j$ satisfies $R_j\le1-h_j$.
Lemma~\ref{lem:gen-balancing} supplies parameters for which
\[
 X_j(r)\le R_j\le1-h_j\quad\text{and}
 \quad X_k(r)\le1\quad \text{for all } k\ne j.
\]
Lemma~\ref{lem:uncovered-point} contradicts the covering.
\end{proof}

\begin{lemma}\label{lem:gen-balancing}
Suppose that $\sum_iw_i=1$ and every $c_{ik}$ is positive.
For each coordinate $j$, there are parameters with $r_j=1$
and $r_k\ge1$ for $k\ne j$ such that
\[
 X_j(r)\le R_j\quad \text{and}
 \quad X_k(r)\le1 \quad\text{for all }k\ne j.
\]
\end{lemma}

\begin{proof}
Consider the set of all parameter vectors $r$ satisfying
\[
 r_j=1,\qquad r_k\ge1,
 \qquad X_j(r)\ge R_j,\quad \text{and}\quad \qquad X_k(r)\le1\quad\text{for all }k\ne j.
\]
This set contains $(1,\ldots,1)$, since
\begin{align*}    
 X_k(1,\ldots,1)=&\sum_iw_i c_{ik}<\sum_iw_i=1
 \quad\text{for all }k,
\\
X_j(1,\ldots,1)=&\sum_iw_i c_{ij}
 \ge\sum_i\frac{w_i c_{ij}}{\sum_{\ell=1}^d c_{i\ell}}=R_j.
\end{align*}
The set is closed and bounded: its defining inequalities imply
\[
 r_kR_j\le r_kX_j(r)
 =\sum_i\frac{w_i r_k c_{ij}}{\|A_r n_i\|^2}
 \le\sum_i\frac{w_i r_k c_{ij}}{r_kc_{ik}^2}
 =\sum_i\frac{w_i c_{ij}}{c_{ik}^2}
 \quad\text{for all }k\ne j,
\]
therefore
\[
 1\le r_k\le\frac1{R_j}\sum_i\frac{w_i c_{ij}}{c_{ik}^2}
 \quad\text{for all }k\ne j.
\]

Choose $r$ in this compact set minimizing $X_j(r)$.
If $X_j(r)=R_j$, the conclusion follows.
Otherwise $X_j(r)>R_j$.
If $X_k(r)<1$ for some $k\ne j$, a sufficiently small increase
of $r_k$ keeps $X_j(r)>R_j$ and $X_k(r)<1$,
while decreasing every span other than $X_k(r)$.
This contradicts minimality, so $X_k(r)=1$ for every $k\ne j$.

Group coordinate $j$ against the remaining coordinates, setting
\[
 a_i:=c_{ij},\qquad b_i:=\sum_{k\ne j}c_{ik}.
\]
Then $a_i,b_i>0$ and
\[
 a_i^2+b_i^2
 =c_{ij}^2+\Bigl(\sum_{k\ne j}c_{ik}\Bigr)^2
 \ge c_{ij}^2+\sum_{k\ne j}c_{ik}^2=1.
\]
The Cauchy--Schwarz inequality gives
\[
 \|A_rn_i\|^2=c_{ij}^2+\sum_{k\ne j}r_kc_{ik}^2
 \ge c_{ij}^2+
 \frac{(\sum_{k\ne j}c_{ik})^2}{\sum_{k\ne j}r_k^{-1}}
 =a_i^2+\lambda b_i^2,
 \quad
 \text{where }
 \lambda:=\left(\sum_{k\ne j}r_k^{-1}\right)^{-1}.
\]
The coefficients $\lambda/r_k$, $k\ne j$, sum to one.
Thus the equalities $X_k(r)=1$ give
\begin{equation}\label{eq:coordinate-balance}
 \sum_iw_i\frac{\lambda b_i}{\|A_rn_i\|^2}
 =\sum_{k\ne j}\frac{\lambda}{r_k}X_k(r)=1.
\end{equation}
Since every $r_k\ge1$, we have $\|A_rn_i\|^2\ge1$.
Moreover, $\|A_rn_i\|^2\ge a_i^2+\lambda b_i^2>\lambda b_i^2$, and hence
\[
 1=\sum_iw_i\frac{\lambda b_i}{\|A_rn_i\|^2}
 <\sqrt{\lambda}\sum_i\frac{w_i}{\|A_rn_i\|}
 \le\sqrt{\lambda}.
\]
Thus $\lambda>1$.
Using $\|A_rn_i\|^2\ge a_i^2+\lambda b_i^2$, we obtain
\[
 \frac{(\lambda-1)a_i+\lambda b_i}{\|A_rn_i\|^2}
 \le\frac{(\lambda-1)a_i+\lambda b_i}{a_i^2+\lambda b_i^2}
 \overset{\eqref{eq:two-variable-bound}}{\le}
1+(\lambda-1)\frac{a_i}{a_i+b_i}.
\]
Multiplying by $w_i$ and summing, using $\sum_iw_i=1$, we obtain a contradiction
with $X_j(r)>R_j$:
\[
(\lambda-1)X_j(r)+1
 \overset{\eqref{eq:coordinate-balance}}{=}
 (\lambda-1)X_j(r)
 +\sum_iw_i\frac{\lambda b_i}{\|A_rn_i\|^2}
 \le1+(\lambda-1)R_j.
  \qedhere
\]
\end{proof}

\subsection*{Funding}
E.B. is supported by the DNRF Chair grant 20231101-28697.
A.P. is supported by the NSF grant DMS 2349045.

\subsection*{Declaration on the Use of AI}
The idea of applying Bang's lemma after a suitable affine change
of coordinates was proposed by the authors and implemented
by ChatGPT.

\bibliographystyle{alpha}
\bibliography{planks}

@article{ball1991plank,
  author  = {Ball, Keith},
  title   = {The plank problem for symmetric bodies},
  journal = {Inventiones Mathematicae},
  volume  = {104},
  year    = {1991},
  pages   = {535--543},
  doi     = {10.1007/BF01245089}
}

@article{bang1950covering,
  author  = {Bang, Th{\o}ger},
  title   = {On covering by parallel-strips},
  journal = {Matematisk Tidsskrift. B},
  volume  = {1950},
  year    = {1950},
  pages   = {49--53}
}

@article{bang1951solution,
  author  = {Bang, Th{\o}ger},
  title   = {A solution of the {``}plank problem{''}},
  journal = {Proceedings of the American Mathematical Society},
  volume  = {2},
  year    = {1951},
  pages   = {990--993}
}

@incollection{bezdek2013tarski,
  author       = {Bezdek, K{\'a}roly},
  title        = {Tarski's plank problem revisited},
  booktitle    = {Geometry---Intuitive, Discrete, and Convex},
  editor       = {B{\'a}r{\'a}ny, Imre and B{\"o}r{\"o}czky, K{\'a}roly J. and Fejes T{\'o}th, G{\'a}bor and Pach, J{\'a}nos},
  series       = {Bolyai Society Mathematical Studies},
  volume       = {24},
  pages        = {45--64},
  publisher    = {Springer},
  address      = {Berlin, Heidelberg},
  year         = {2013},
  doi          = {10.1007/978-3-642-41498-5_2},
  eprint       = {0903.4637},
  archivePrefix= {arXiv},
  primaryClass = {math.MG}
}

@article{tarski1932uwagi,
  author  = {Tarski, Alfred},
  title   = {Uwagi o stopniu r{\'o}wnowa{\.z}no{\'s}ci wielok{\k{a}}t{\'o}w},
  journal = {Parametr},
  volume  = {2},
  year    = {1932},
  pages   = {310--314},
  note    = {In Polish}
}

@article{bezdek2003covering,
  title={Covering an annulus by strips},
  author={Bezdek, Andr{\'a}s},
  journal={Discrete \& Computational Geometry},
  volume={30},
  number={2},
  pages={177--180},
  year={2003},
  publisher={Springer},
    note    = {U.S.-Hungarian Workshops on Discrete Geometry and Convexity
             (Budapest, 1999/Auburn, AL, 2000)}
}

@article{white2007covering,
  title={Covering polygonal annuli by strips},
  author={White, Stuart and Wisewell, Laura},
  journal={Discrete \& Computational Geometry},
  volume={37},
  number={4},
  pages={577--585},
  year={2007},
  publisher={Springer}
}

@article{smurov2010covering,
  title={Covering a cube with holes by strips},
  author={Smurov, M.V. and Bogataya, S.I. and Bogatyi, S.A.},
  journal={Topology and its Applications},
  volume={157},
  number={4},
  pages={708--710},
  year={2010},
  publisher={Elsevier}
}

@article{ambrus2025covering,
  title={Covering Spiky Annuli by Planks},
  author={Ambrus, Gergely and Huddell, Julian and Lai, Maggie and Quirk, Matthew and Williams, Elias},
  journal={Discrete \& Computational Geometry},
  volume  = {76},
  number  = {2},
  pages   = {978--987},
  year    = {2026},
  publisher={Springer}
}

@book{brass2005research,
  title={Research problems in discrete geometry},
  author={Brass, Peter and Moser, William OJ and Pach, J{\'a}nos},
  year={2005},
  publisher={Springer}
}

@misc{ambrus2026plankcoveringscubeannuli,
      title={On plank coverings of cube annuli}, 
      author={Gergely Ambrus and Matthew Corbo and Dakota Nobuyuki Goto and Desmond Leitz and Kai Mawhinney},
      year={2026},
      note={arXiv:2610.18770},
      archivePrefix={arXiv},
      primaryClass={math.MG},
      url={https://arxiv.org/abs/2610.18770}, 
}

@misc{bakaev2026volumestabilityplankcovers,
      title={Volume stability of plank covers}, 
      author={Egor Bakaev and Amir Yehudayoff},
      year={2026},
      note={arXiv:2609.30409},
      archivePrefix={arXiv},
      primaryClass={math.MG},
      url={https://arxiv.org/abs/2609.30409}, 
}

@misc{akopyan2026shortproofballsplank,
      title={A short proof of {Ball}'s plank theorem}, 
      author={Arseniy Akopyan and Alexander Polyanskii},
      year={2026},
      note={arXiv:2609.36163},
      archivePrefix={arXiv},
      primaryClass={math.MG},
      url={https://arxiv.org/abs/2609.36163}, 
}

\end{document}